\documentclass[12pt]{article}

\usepackage{amsmath,amssymb,amsthm,amscd,a4wide}

\begin{document}

\newcommand{\ad}{{\rm ad}}
\newcommand{\cri}{{\rm cri}}
\newcommand{\End}{{\rm{End}\ts}}
\newcommand{\Rep}{{\rm{Rep}\ts}}
\newcommand{\Hom}{{\rm{Hom}}}
\newcommand{\Mat}{{\rm{Mat}}}
\newcommand{\ch}{{\rm{ch}\ts}}
\newcommand{\chara}{{\rm{char}\ts}}
\newcommand{\diag}{{\rm diag}}
\newcommand{\non}{\nonumber}
\newcommand{\wt}{\widetilde}
\newcommand{\wh}{\widehat}
\newcommand{\ot}{\otimes}
\newcommand{\la}{\lambda}
\newcommand{\La}{\Lambda}
\newcommand{\De}{\Delta}
\newcommand{\al}{\alpha}
\newcommand{\be}{\beta}
\newcommand{\ga}{\gamma}
\newcommand{\Ga}{\Gamma}
\newcommand{\ep}{\epsilon}
\newcommand{\ka}{\kappa}
\newcommand{\vk}{\varkappa}
\newcommand{\si}{\sigma}
\newcommand{\vs}{\varsigma}
\newcommand{\vp}{\varphi}
\newcommand{\de}{\delta}
\newcommand{\ze}{\zeta}
\newcommand{\om}{\omega}
\newcommand{\ee}{\epsilon^{}}
\newcommand{\su}{s^{}}
\newcommand{\hra}{\hookrightarrow}
\newcommand{\ve}{\varepsilon}
\newcommand{\ts}{\,}
\newcommand{\vac}{\mathbf{1}}
\newcommand{\di}{\partial}
\newcommand{\qin}{q^{-1}}
\newcommand{\tss}{\hspace{1pt}}
\newcommand{\Sr}{ {\rm S}}
\newcommand{\U}{ {\rm U}}
\newcommand{\BL}{ {\overline L}}
\newcommand{\BE}{ {\overline E}}
\newcommand{\BP}{ {\overline P}}
\newcommand{\AAb}{\mathbb{A}\tss}
\newcommand{\CC}{\mathbb{C}\tss}
\newcommand{\KK}{\mathbb{K}\tss}
\newcommand{\QQ}{\mathbb{Q}\tss}
\newcommand{\SSb}{\mathbb{S}\tss}
\newcommand{\ZZ}{\mathbb{Z}\tss}
\newcommand{\X}{ {\rm X}}
\newcommand{\Y}{ {\rm Y}}
\newcommand{\Z}{{\rm Z}}
\newcommand{\Ac}{\mathcal{A}}
\newcommand{\Lc}{\mathcal{L}}
\newcommand{\Mc}{\mathcal{M}}
\newcommand{\Pc}{\mathcal{P}}
\newcommand{\Qc}{\mathcal{Q}}
\newcommand{\Rc}{\mathcal{R}}
\newcommand{\Sc}{\mathcal{S}}
\newcommand{\Tc}{\mathcal{T}}
\newcommand{\Bc}{\mathcal{B}}
\newcommand{\Ec}{\mathcal{E}}
\newcommand{\Fc}{\mathcal{F}}
\newcommand{\Hc}{\mathcal{H}}
\newcommand{\Uc}{\mathcal{U}}
\newcommand{\Vc}{\mathcal{V}}
\newcommand{\Wc}{\mathcal{W}}
\newcommand{\Yc}{\mathcal{Y}}
\newcommand{\Ar}{{\rm A}}
\newcommand{\Br}{{\rm B}}
\newcommand{\Ir}{{\rm I}}
\newcommand{\Fr}{{\rm F}}
\newcommand{\Jr}{{\rm J}}
\newcommand{\Or}{{\rm O}}
\newcommand{\GL}{{\rm GL}}
\newcommand{\Spr}{{\rm Sp}}
\newcommand{\Rr}{{\rm R}}
\newcommand{\Zr}{{\rm Z}}
\newcommand{\gl}{\mathfrak{gl}}
\newcommand{\middd}{{\rm mid}}
\newcommand{\ev}{{\rm ev}}
\newcommand{\Pf}{{\rm Pf}}
\newcommand{\Norm}{{\rm Norm\tss}}
\newcommand{\oa}{\mathfrak{o}}
\newcommand{\spa}{\mathfrak{sp}}
\newcommand{\osp}{\mathfrak{osp}}
\newcommand{\f}{\mathfrak{f}}
\newcommand{\g}{\mathfrak{g}}
\newcommand{\h}{\mathfrak h}
\newcommand{\n}{\mathfrak n}
\newcommand{\z}{\mathfrak{z}}
\newcommand{\Zgot}{\mathfrak{Z}}
\newcommand{\p}{\mathfrak{p}}
\newcommand{\sll}{\mathfrak{sl}}
\newcommand{\agot}{\mathfrak{a}}
\newcommand{\qdet}{ {\rm qdet}\ts}
\newcommand{\Ber}{ {\rm Ber}\ts}
\newcommand{\HC}{ {\mathcal HC}}
\newcommand{\cdet}{{\rm cdet}}
\newcommand{\rdet}{{\rm rdet}}
\newcommand{\tr}{ {\rm tr}}
\newcommand{\gr}{\mathbf{gr}}
\newcommand{\str}{ {\rm str}}
\newcommand{\loc}{{\rm loc}}
\newcommand{\Gr}{{\rm G}}
\newcommand{\sgn}{ {\rm sgn}\ts}
\newcommand{\sign}{{\rm sgn}}
\newcommand{\ba}{\bar{a}}
\newcommand{\bb}{\bar{b}}
\newcommand{\bi}{\bar{\imath}}
\newcommand{\bj}{\bar{\jmath}}
\newcommand{\bk}{\bar{k}}
\newcommand{\bl}{\bar{l}}
\newcommand{\hb}{\mathbf{h}}
\newcommand{\Sym}{\mathfrak S}
\newcommand{\fand}{\quad\text{and}\quad}
\newcommand{\Fand}{\qquad\text{and}\qquad}
\newcommand{\For}{\qquad\text{or}\qquad}
\newcommand{\OR}{\qquad\text{or}\qquad}

\renewcommand{\theequation}{\arabic{section}.\arabic{equation}}

\newtheorem{thm}{Theorem}[section]
\newtheorem{lem}[thm]{Lemma}
\newtheorem{prop}[thm]{Proposition}
\newtheorem{cor}[thm]{Corollary}
\newtheorem{conj}[thm]{Conjecture}
\newtheorem*{mthm}{Main Theorem}
\newtheorem*{mthma}{Theorem A}
\newtheorem*{mthmb}{Theorem B}

\theoremstyle{definition}
\newtheorem{defin}[thm]{Definition}

\theoremstyle{remark}
\newtheorem{remark}[thm]{Remark}
\newtheorem{example}[thm]{Example}

\newcommand{\bth}{\begin{thm}}
\renewcommand{\eth}{\end{thm}}
\newcommand{\bpr}{\begin{prop}}
\newcommand{\epr}{\end{prop}}
\newcommand{\ble}{\begin{lem}}
\newcommand{\ele}{\end{lem}}
\newcommand{\bco}{\begin{cor}}
\newcommand{\eco}{\end{cor}}
\newcommand{\bde}{\begin{defin}}
\newcommand{\ede}{\end{defin}}
\newcommand{\bex}{\begin{example}}
\newcommand{\eex}{\end{example}}
\newcommand{\bre}{\begin{remark}}
\newcommand{\ere}{\end{remark}}
\newcommand{\bcj}{\begin{conj}}
\newcommand{\ecj}{\end{conj}}

\newcommand{\bal}{\begin{aligned}}
\newcommand{\eal}{\end{aligned}}
\newcommand{\beq}{\begin{equation}}
\newcommand{\eeq}{\end{equation}}
\newcommand{\ben}{\begin{equation*}}
\newcommand{\een}{\end{equation*}}

\newcommand{\bpf}{\begin{proof}}
\newcommand{\epf}{\end{proof}}

\def\beql#1{\begin{equation}\label{#1}}

\title{\Large\bf Eigenvalues of Bethe vectors in the Gaudin model}

\author{{A. I. Molev\quad and\quad E. E. Mukhin}}

\date{} 
\maketitle

\vspace{30 mm}

\begin{abstract}
A theorem of Feigin, Frenkel and Reshetikhin provides expressions
for the eigenvalues of the higher Gaudin Hamiltonians on the Bethe vectors
in terms of elements of the center of the affine vertex algebra at the critical level.
In our recent work, explicit Harish-Chandra images of generators of the center
were calculated in all classical types. We combine these results to calculate
the eigenvalues of the higher Gaudin Hamiltonians on
the Bethe vectors in an explicit form.
The Harish-Chandra images can be interpreted as elements of classical $\Wc$-algebras.
We provide a direct connection between the rings of $q$-characters and
classical $\Wc$-algebras by calculating classical limits of the corresponding
screening operators.

\end{abstract}


\vspace{35 mm}

\noindent
School of Mathematics and Statistics\newline
University of Sydney,
NSW 2006, Australia\newline
alexander.molev@sydney.edu.au

\vspace{7 mm}

\noindent
Department of Mathematical Sciences\newline
Indiana University -- Purdue University Indianapolis\newline
402 North Blackford St, Indianapolis, IN 46202-3216, USA\newline
mukhin@math.iupui.edu

\newpage

%

\section{Introduction}
\label{sec:int}
\setcounter{equation}{0}

In their seminal paper \cite{ffr:gm}, Feigin, Frenkel and Reshetikhin
established a connection between the center $\z(\wh\g)$ of the affine vertex algebra
at the critical level and the higher Gaudin Hamiltonians. They used the
Wakimoto modules over the affine Kac--Moody algebra $\wh\g$ to calculate
the eigenvalues of the Hamiltonians on the Bethe vectors of the Gaudin model
associated with an arbitrary simple Lie algebra $\g$. The calculation depends on
the choice of an element $S$ of the center and the result is written
in terms of the Harish-Chandra image of $S$; see also \cite{ffr:oi},
\cite{fft:gm} and \cite{r:si}
for a relationship with the opers and generalizations to non-homogeneous
Hamiltonians.

The center $\z(\wh\g)$ is a commutative associative algebra whose structure was
described by a theorem of Feigin and Frenkel~\cite{ff:ak}, which states that
$\z(\wh\g)$ is an algebra of polynomials in infinitely many variables; see \cite{f:lc} for
a detailed exposition of these results. Simple explicit formulas for generators
of this algebra were found in \cite{ct:qs} for type $A$ and in \cite{m:ff}
for types $B$, $C$ and $D$; see also \cite{cm:ho} and \cite{mr:mm} for simpler arguments
in type $A$ and extensions to Lie superalgebras.
The calculation of the Harish-Chandra images of the generators in type $A$
is straightforward, whereas types $B$, $C$ and $D$ require a rather involved
application of the $q$-characters; see \cite{mm:yc}.
Our goal in this paper is to apply these results
to get the action of the higher Gaudin Hamiltonians on tensor products
of representations of $\g$ in an explicit form and calculate the corresponding
eigenvalues of the Bethe vectors. In type $A$ we thus reproduce the results of
\cite{mtv:be} obtained by a different method based on the Bethe ansatz.

We will begin with a brief exposition of some results of \cite{ffr:gm} and \cite{fft:gm}.
Our main focus will be on Theorem~6.7 from \cite{fft:gm} expressing eigenvalues
of a generalized Gaudin algebra on Bethe vectors in terms of opers associated
with tensor products of Verma modules. Then we will apply this theorem
to the classical Lie algebras to write explicit Gaudin operators and their
eigenvalues on Bethe vectors.

A connection between the Yangian characters (or $q$-characters) and
the Segal--Sugawara operators played an essential role in the calculation
of the Harish-Chandra images in \cite{mm:yc}. We will explore this connection further
by constructing a map
$\gr$ taking a character to an element of the associated classical
$\Wc$-algebra. We will also
establish multiplicativity and surjectivity properties of this map.

\section{Feigin--Frenkel center and Bethe vectors}
\label{sec:ga}
\setcounter{equation}{0}

Let $\g$ be a simple Lie algebra over $\CC$ equipped with
a standard symmetric invariant bilinear form $\langle\ts\ts,\ts\rangle$
defined as a normalized Killing form
\beql{killi}
\langle X,Y\rangle=\frac{1}{2\tss h^{\vee}}\ts\tr\ts\big(\ad\tss X\ts\ad\tss Y\big),
\eeq
where $h^{\vee}$ is the {\it dual Coxeter number\/} for $\g$.
The corresponding {\it affine Kac--Moody algebra\/} $\wh\g$
is defined as the central
extension
\beql{km}
\wh\g=\g\tss[t,t^{-1}]\oplus\CC K,
\eeq
where $\g[t,t^{-1}]$ is the Lie algebra of Laurent
polynomials in $t$ with coefficients in $\g$; see \cite{k:id}.
For any $r\in\ZZ$ and $X\in\g$
we set $X[r]=X\ts t^r$. The commutation relations of the Lie algebra $\wh\g$
have the form
\ben
\big[X[r],Y[s]\big]=[X,Y][r+s]+r\ts\de_{r,-s}\langle X,Y\rangle\ts K,
\qquad X, Y\in\g,
\een
and the element $K$ is central in $\wh\g$.

The universal enveloping algebra {\em at the critical level}
$\U_{\cri}(\wh\g)$ is the quotient of $\U(\wh\g)$
by the ideal generated by $K+h^{\vee}$.
Let $\Ir$ denote the left ideal of $\U_{\cri}(\wh\g)$ generated by $\g[t]$
and let $\Norm\tss\Ir$ be its normalizer,
\ben
\Norm\tss\Ir=\{v\in \U_{\cri}(\wh\g)\ |\ \Ir\tss v\subseteq \Ir\}.
\een
The normalizer is a subalgebra of $\U_{\cri}(\wh\g)$, and $\Ir$
is a two-sided ideal of $\Norm\tss\Ir$.
The {\it Feigin--Frenkel center\/} $\z(\wh\g)$ is the associative algebra
defined as the quotient
\beql{ffnorm}
\z(\wh\g)=\Norm\tss\Ir/\Ir.
\eeq
Any element of $\z(\wh\g)$
is called a {\it Segal--Sugawara vector\/}.
The quotient
\beql{vacm}
V_{\cri}(\g)=\U_{\cri}(\wh\g)/\Ir
\eeq
is the {\em vacuum module} at the critical level over $\wh\g$.
It possesses a vertex algebra structure. As a vector space,
$V_{\cri}(\g)$ is isomorphic to the universal enveloping algebra
$\U(\wh\g_-)$, where $\wh\g_-=t^{-1}\g[t^{-1}]$.
Hence, we have
a vector space embedding
\ben
\z(\wh\g)\hra \U(\wh\g_-).
\een
Since $\U(\wh\g_-)$ is a subalgebra of $\U_{\cri}(\wh\g)$,
the embedding
is an algebra homomorphism so that the Feigin--Frenkel center
$\z(\wh\g)$ can be regarded as a
subalgebra of $\U(\wh\g_-)$. This subalgebra is commutative
which is seen from its identification with the {\em center}
of the vertex algebra $V_{\cri}(\g)$ by
\beql{ffvac}
\z(\wh\g)=\{v\in V_{\cri}(\g)\ |\ \g[t]\tss v=0\}.
\eeq
As a vertex algebra, the vacuum module $V_{\cri}(\g)$ is equipped
with the translation operator
\beql{transl}
T:V_{\cri}(\g)\to V_{\cri}(\g),
\eeq
which is determined by the properties
\ben
T:1\mapsto 0\Fand
\big[T,X[r]\big]=-r\tss X[r - 1],\quad X\in\g.
\een
We also regard $T$ as a derivation of the algebra $\U(\wh\g_-)$.
Its subalgebra $\z(\wh\g)$ is $T$-invariant.
By the Feigin--Frenkel theorem, there exist
elements $S_1,\dots,S_n\in \z(\wh\g)$, where $n={\rm rank\ts}\g$, such that
all elements $T^{\tss r}S_l$ are algebraically
independent, and every Segal--Sugawara vector is a polynomial
in these elements:
\beql{genz}
\z(\wh\g)=\CC[T^{\tss r}S_l\ |\ l=1,\dots,n,\ \ r\geqslant 0].
\eeq
We call such a family $S_1,\dots,S_n$ a {\em complete set of Segal--Sugawara vectors}.

Choose a Cartan subalgebra $\h$ of $\g$ and
a triangular decomposition $\g=\n_-\oplus\h\oplus \n_+$.
Consider $\U(\wh\g_-)$ as the adjoint $\g$-module by
regarding $\g$ as the span of the elements $X[0]$ with $X\in\g$.
Denote by $\U(\wh\g_-)^{\h}$ the
subalgebra of $\h$-invariants under this action.
Consider the left ideal $\Jr$ of the algebra $\U(\wh\g_-)$ generated by
all elements $X[r]$ with $X\in\n_-$ and $r< 0$.
By the Poincar\'e--Birkhoff--Witt
theorem, the intersection
$\U(\wh\g_-)^{\h}\cap \Jr$ is a two-sided ideal of
$\U(\wh\g_-)^{\h}$ and we have a direct sum decomposition
\ben
\U(\wh\g_-)^{\h}=\big(\U(\wh\g_-)^{\h}\cap \Jr\big)\oplus\U(\wh\h_-),
\een
where $\wh\h_-=t^{-1}\h[t^{-1}]$. The projection to the second summand
is a homomorphism
\beql{hchaff}
\U(\wh\g_-)^{\h}\to \U(\wh\h_-)
\eeq
which is an affine version of the {\em Harish-Chandra homomorphism}.
By the Feigin--Frenkel theorem, the restriction
of the homomorphism \eqref{hchaff} to
the subalgebra $\z(\wh\g)$ yields an isomorphism
\beql{hchiaff}
\f:\z(\wh\g)\to \Wc({}^L\g),
\eeq
where $\Wc({}^L\g)$ is the {\it classical $\Wc$-algebra\/} associated with the
Langlands dual Lie algebra ${}^L\g$; see \cite{f:lc} for
a detailed exposition of these results.
The $\Wc$-algebra $\Wc({}^L\g)$ can be defined as a subalgebra
of $\U(\wh\h_-)$ which consists of the elements
annihilated by the {\it screening operators\/}; see \cite[Sec.~8.1]{f:lc}
and also \cite{mm:yc} for explicit formulas in the classical types.

Given any element $\chi\in\g^*$ and
a nonzero $z\in\CC$, the mapping
\beql{evalr}
\U(\wh\g_-)\to \U(\g),
\qquad X[r]\mapsto X\tss z^r+\de_{r,-1}\ts\chi(X),\quad X\in\g,
\eeq
defines an algebra homomorphism.
Using the coassociativity of the standard coproduct on $\U(\wh\g_-)$
defined by
\ben
\Delta: Y\mapsto Y\ot 1+1\ot Y,\qquad Y\in\wh\g_-,
\een
for any $\ell\geqslant 1$ we get the homomorphism
\beql{comult}
\U(\wh\g_-)\to \U(\wh\g_-)^{\ot\tss\ell}
\eeq
as an iterated coproduct map. Now fix distinct complex numbers
$z_1,\dots,z_\ell$ and let $u$ be a complex parameter.
Applying homomorphisms of the form \eqref{evalr} to the tensor factors in
\eqref{comult}, we get another homomorphism
\beql{psiu}
\Psi:\U(\wh\g_-)\to \U(\g)^{\ot\tss\ell},
\eeq
given by
\ben
\Psi:X[r]\mapsto \sum_{a=1}^\ell X_a(z_a-u)^r+\de_{r,-1}\ts\chi(X)\in \U(\g)^{\ot\tss\ell},
\een
where $X_a=1^{\ot (a-1)}\ot X\ot 1^{\ot (\ell-a)}$; see \cite{r:si}.
We will twist this homomorphism by the involutive anti-automorphism
\beql{sgn}
\vs: \U(\wh\g_-)\to \U(\wh\g_-),
\qquad X[r]\mapsto -X[r],\quad X\in\g,
\eeq
to get the anti-homomorphism
\beql{phiu}
\Phi:\U(\wh\g_-)\to \U(\g)^{\ot\tss\ell},
\eeq
defined as the composition $\Phi=\Psi\circ\vs$.
Since $\z(\wh\g)$ is a commutative subalgebra of $\U(\wh\g_-)$,
the image of $\z(\wh\g)$ under $\Phi$ is a commutative subalgebra
$\Ac(\g)_{\chi}$ of $\U(\g)^{\ot\tss\ell}$, depending on the chosen
parameters $z_1,\dots,z_\ell$, but it does not depend on $u$
\cite{r:si}; see also \cite[Sec.~2]{fft:gm}.

Introduce the
standard Chevalley generators $e_i,h_i,f_i$ with $i=1,\dots,n$
of the simple Lie algebra $\g$ of rank $n$. The generators $h_i$ form
a basis of the Cartan subalgebra $\h$ of $\g$, while the $e_i$ and $f_i$
generate the respective nilpotent subalgebras $\n_+$ and $\n_-$.
Let $A=[a_{ij}]$
be the Cartan matrix of $\g$ so that the defining relations of $\g$ take the form
\begin{alignat}{2}
[e_i,f_j]&=\de_{ij}h_i,\qquad [h_i,h_j]&&=0,
\non\\
[h_i,e_j]&=a_{ij}\tss e_j,\qquad [h_i,f_j]&&=-a_{ij}\tss f_j,
\non
\end{alignat}
together with the Serre relations; see e.g. \cite{k:id}.
Given $\la\in\h^*$, the Verma module $M_{\la}$ is the quotient of $\U(\g)$
by the left ideal generated by $\n_+$ and the elements $h_i-\la(h_i)$
with $i=1,\dots,n$. We denote the image of $1$ in $M_{\la}$ by $1_{\la}$.

For any weights $\la_1,\dots,\la_\ell\in\h^*$ consider the tensor product
of the Verma modules $M_{\la_1}\ot\dots\ot M_{\la_\ell}$. We will now describe common
eigenvectors for the commutative subalgebra $\Ac(\g)_{\chi}$
in this tensor product. For a set of distinct complex numbers $w_1,\dots, w_m$
with $w_i\ne z_j$ and a collection (multiset) of labels $i_1,\dots,i_m\in\{1,\dots,n\}$
introduce the {\em Bethe vector}
\ben
\phi(w_1^{i_1},\dots,w_m^{i_m})\in M_{\la_1}\ot\dots\ot M_{\la_\ell}
\een
by the following formula which originates in \cite{sv:ah}; see \cite{bf:os} and
also \cite{ffr:gm}, \cite{mv:nb} and references therein:
\beql{bv}
\phi(w_1^{i_1},\dots,w_m^{i_m})=\sum_{(I^1,\dots,I^\ell)}
\bigotimes_{k=1}^\ell\ts\prod_{s=1}^{a_k}\ts\frac{1}{w_{j^k_s}-w_{j^k_{s+1}}}
\ts\prod_{r\in I^k}f_{i_r}\ts 1_{\la_k},
\eeq
summed over all ordered partitions $I^1\cup I^2\cup \dots \cup I^\ell$ of the set $\{1,\dots,m\}$
into ordered subsets $I^k=\{j^k_1,j^k_2,\dots,j^k_{a_k}\}$ with the products taken from left to right,
where $w_{j^k_{s+1}}:=z_k$ for $s=a_k$.

Now suppose that $\chi\in\h^*$. We regard $\chi$ as a functional on
$\g$ which vanishes on $\n_+$ and $\n_-$. The system
of the {\em Bethe ansatz equations} takes the form
\beql{bae}
\sum_{i=1}^\ell\frac{\langle \check\al_{i_j},\la_i\rangle}{w_j-z_i}-\sum_{s\ne j}
\frac{\langle \check\al_{i_j},\al_{i_s}\rangle}{w_j-w_s}=\langle \check\al_{i_j},\chi\rangle,
\qquad j=1,\dots,m,
\eeq
where the $\al_l$ and $\check\al_l$ denote the simple roots and coroots, respectively; see \cite{k:id}.

We are now in a position to describe the eigenvalues of the Gaudin Hamiltonians
on the Bethe vectors. Given the above parameters,
introduce the homomorphism
from $\U(\wh\h_-)$ to rational functions in $u$ by the rule:
\beql{hw}
\varrho:H[-r-1]\mapsto \frac{\di_u^{\ts r}}{r!}\ts \Hc(u),\qquad H\in\h,\quad
r\geqslant 0,
\eeq
where
\ben
\Hc(u)=\sum_{a=1}^\ell \frac{\la_a(H)}
{u-z_a}-\sum_{j=1}^m\frac{\al_{i_j}(H)}{u-w_j}-\chi(H).
\een
Let $S\in \z(\wh\g)$ be a Segal--Sugawara vector. The composition $\varrho\circ\f$
of this homomorphism with the isomorphism \eqref{hchiaff} takes $S$ to
a rational function $\varrho\big(\f(S)\big)$ in $u$.
Furthermore,
we regard the image $\Phi(S)$ of $S$ under the anti-homomorphism \eqref{phiu}
as an operator in the tensor product of Verma modules
$M_{\la_1}\ot\dots\ot M_{\la_\ell}$.
The following is essentially a reformulation of Theorems~6.5 and 6.7 from \cite{fft:gm};
in the case $\chi=0$ the result goes back to \cite[Theorem~3]{ffr:gm}.

\bth\label{thm:eigen}
Suppose that the Bethe ansatz equations \eqref{bae} are satisfied.
If the Bethe vector $\phi(w_1^{i_1},\dots,w_m^{i_m})$ is nonzero, then it is an eigenvector
for the operator $\Phi(S)$ with the eigenvalue $\varrho\big(\f(S)\big)$.
\qed
\eth

In what follows we will rely on the results of \cite{cm:ho}, \cite{m:ff} and \cite{mm:yc}
to give explicit formulas for the operators $\Phi(S_i)$
and their eigenvalues $\varrho\big(\f(S_i)\big)$ on the Bethe vectors for
complete sets of Segal--Sugawara vectors $S_1,\dots,S_n$
in all classical types.

\section{Gaudin Hamiltonians and eigenvalues}
\label{sec:gh}
\setcounter{equation}{0}

We will use the extended Lie algebra $\wh\g\oplus\CC\tau$ where the element $\tau$
satisfies the commutation relations
\beql{taur}
\big[\tau,X[r]\tss\big]=-r\ts X[r-1],\qquad
\big[\tau,K\big]=0.
\eeq
Consider the extension of \eqref{hchiaff}
to the isomorphism
\beql{hchiaffe}
\f:\z(\wh\g)\ot\CC[\tau]\to \Wc({}^L\g)\ot\CC[\tau],
\eeq
which is identical on $\CC[\tau]$.

For an arbitrary $N\times N$
matrix $A=[A_{ij}]$ with entries in a ring
we define its {\em column-determinant} $\cdet\ts A$
and {\em row-determinant} $\rdet\ts A$
by the respective formulas
\beql{cdetdef}
\cdet\ts A=\sum_{\si\in\Sym_N}\sgn\si\cdot A_{\si(1)1}\dots
A_{\si(N)N}
\eeq
and
\beql{rdetdef}
\rdet\ts A=\sum_{\si\in\Sym_N}\sgn\si\cdot A_{1\si(1)}\dots
A_{N\si(N)},
\eeq
where $\Sym_N$ denotes the symmetric group.

\subsection{Type $A$}
\label{subssec:typea}

We will work with the reductive Lie algebra $\gl_N$ rather than the simple
Lie algebra $\sll_N$ of type $A$. We let $E_{ij}$ with $i,j=1,\dots,N$ be
the standard basis of $\gl_N$.
Denote by $\h$, $\n_+$ and $\n_-$
the subalgebras of $\gl_N$ spanned by the diagonal, upper-triangular
and lower-triangular matrices, respectively, so that $E_{11},\dots,E_{NN}$
is a basis of $\h$.

We start by recalling the constructions of some complete sets of Segal--Sugawara vectors
for $\gl_N$.
For each $a\in\{1,\dots,m\}$
introduce the element $E[r]_a$ of the algebra
\beql{tenprka}
\underbrace{\End\CC^{N}\ot\dots\ot\End\CC^{N}}_m{}\ot\U
\eeq
by
\beql{matnota}
E[r]_a=\sum_{i,j=1}^{N}
1^{\ot(a-1)}\ot e_{ij}\ot 1^{\ot(m-a)}\ot E_{ij}[r],
\eeq
where the $e_{ij}$ are the standard matrix units and $\U$ stands
for the universal enveloping algebra of
$\wh\gl_N\oplus\CC\tau$.
Let $H^{(m)}$ and $A^{(m)}$ denote the respective images of the
normalized symmetrizer and anti-symmetrizer in the group algebra
for the symmetric group $\Sym_m$ under
its natural action on $(\CC^{N})^{\ot m}$. In particular,
$H^{(m)}$ and $A^{(m)}$ are idempotents and we identify them with the
respective elements
$H^{(m)}\ot 1$ and $A^{(m)}\ot 1$ of the algebra \eqref{tenprka}.
Define the elements
$\vp^{}_{m\tss a},\psi^{}_{m\tss a}, \theta^{}_{m\tss a}\in
\U\big(t^{-1}\gl_N[t^{-1}]\big)$
by the expansions
\begin{align}\label{deftra}
\tr\ts A^{(m)} \big(\tau+E[-1]_1\big)\dots \big(\tau+E[-1]_m\big)
&=\vp^{}_{m\tss0}\ts\tau^m+\vp^{}_{m\tss1}\ts\tau^{m-1}
+\dots+\vp^{}_{m\tss m},\\[0.7em]
\label{deftrh}
\tr\ts H^{(m)} \big(\tau+E[-1]_1\big)\dots \big(\tau+E[-1]_m\big)
&=\psi^{}_{m\tss0}\ts\tau^m+\psi^{}_{m\tss1}\ts\tau^{m-1}
+\dots+\psi^{}_{m\tss m},
\end{align}
where
the traces are taken with respect to all $m$ copies of $\End\CC^N$
in \eqref{tenprka},
and
\beql{deftracepa}
\tr\ts \big(\tau+E[-1]\big)^m=\theta^{}_{m\tss0}\ts\tau^m+\theta^{}_{m\tss1}\ts\tau^{m-1}
+\dots+\theta^{}_{m\tss m}.
\eeq
Expressions like $\tau+E[-1]$ are understood as matrices, where
$\tau$ is regarded as the scalar matrix $\tau\tss I$.
Furthermore, expand the column-determinant of this matrix as a polynomial in $\tau$,
\beql{coldetal}
\cdet\ts \big(\tau+E[-1]\big)=\tau^N+\vp^{}_{1}\ts\tau^{N-1}
+\dots+\vp^{}_{N},\qquad \vp^{}_{m}\in
\U\big(t^{-1}\gl_N[t^{-1}]\big).
\eeq
We have $\vp^{}_{m\tss m}=\vp^{}_{m}$ for $m=1,\dots,N$.

\bth\label{thm:allff}
All elements $\vp^{}_{m\tss a}$,
$\psi^{}_{m\tss a}$ and $\theta^{}_{m\tss a}$
belong to the Feigin--Frenkel center $\z(\wh\gl_N)$.
Moreover, each of the families
\ben
\vp^{}_{1},\dots,\vp^{}_{N},\qquad \psi^{}_{1\tss 1},\dots,\psi^{}_{N\tss N}
\Fand \theta^{}_{1\tss 1},\dots,\theta^{}_{N\tss N}
\een
is a complete set of Segal--Sugawara vectors for $\gl_N$.
\qed
\eth

This theorem goes back to \cite{ct:qs}, where the elements $\vp^{}_{m}$
were first discovered (in a slightly different form).
A direct proof of the theorem for the coefficients of the polynomial
\eqref{coldetal} was given in \cite{cm:ho}.
The elements $\psi^{}_{m\tss a}$ are related to $\vp^{}_{m\tss a}$
through the quantum MacMahon Master Theorem of \cite{glz:qm}, while
a relationship between the $\vp^{}_{m\tss a}$ and $\theta^{}_{m\tss a}$
is provided by a Newton-type identity given in \cite[Theorem~15]{cfr:ap}.
Note that super-versions of these relations between the families
of Segal--Sugawara vectors for the Lie
superalgebra $\gl_{m|n}$ were given in the paper \cite{mr:mm}, which also provides
simpler arguments in the purely even case.

We will calculate the images of the Segal--Sugawara vectors under
the involution \eqref{sgn}. We extend it to the algebra
$\U\big(t^{-1}\gl_N[t^{-1}]\big)\ot\CC[\tau]$ with the action on $\CC[\tau]$
as the identity map.

\ble\label{lem:sgn}
For the images with respect to the involution $\vs$ we have
\begin{align}
\tr\ts A^{(m)} \big(\tau+E[-1]_1\big)\dots \big(\tau+E[-1]_m\big)
&\mapsto \tr\ts A^{(m)} \big(\tau-E[-1]_1\big)\dots \big(\tau-E[-1]_m\big),
\label{am}\\[0.7em]
\tr\ts H^{(m)} \big(\tau+E[-1]_1\big)\dots \big(\tau+E[-1]_m\big)
&\mapsto \tr\ts H^{(m)} \big(\tau-E[-1]_1\big)\dots \big(\tau-E[-1]_m\big),
\label{hm}\\[0.8em]
\tr\ts \big(\tau+E[-1]\big)^m &\mapsto \tr\ts \big(\tau-E^{\tss t}[-1]\big)^m,
\label{tr}
\end{align}
and
\beql{cdetim}
\cdet\ts \big(\tau+E[-1]\big)\mapsto \cdet\ts \big(\tau-E^{\tss t}[-1]\big),
\eeq
where $t$ denotes the standard matrix transposition.
\ele

\bpf
The left hand side of \eqref{am} equals a linear combination of expressions
of the form
\beql{sumae}
\tr\ts A^{(m)} E[r_1]_{a_1}\dots E[r_p]_{a_p}\ts \tau^k
\eeq
with $1\leqslant a_1<\dots<a_p\leqslant m$. However, such an expression
remains unchanged under any permutation of the factors $E[r_i]_{a_i}$.
This follows from the commutation relations
\ben
E[r]_a\ts E[s]_b-E[s]_b\ts E[r]_a
=P_{a\tss b}\tss E[r+s]_b-E[r+s]_b\tss P_{a\tss b}
\een
for $a<b$, where
\beql{pdef}
P_{a\tss b}=\sum_{i,j=1}^N 1^{\ot(a-1)}\ot e_{ij}
\ot 1^{\ot(b-a-1)}\ot e_{j\tss i}\ot 1^{\ot(m-b)}
\eeq
is the permutation operator.
We only need to observe that $A^{(m)}P_{a\tss b}=P_{a\tss b}A^{(m)}=-A^{(m)}$
and use the cyclic property of trace.
Hence the image of \eqref{sumae}
under $\vs$ equals
\ben
(-1)^p\ts\tr\ts A^{(m)} E[r_1]_{a_1}\dots E[r_p]_{a_p}\ts \tau^k
\een
which verifies \eqref{am}. The same argument proves \eqref{hm}.
Now \eqref{cdetim} follows from the relation
\beql{idtrd}
\cdet\ts \big(\tau+E[-1]\big)=
\tr\ts A^{(N)} \big(\tau+E[-1]_1\big)\dots \big(\tau+E[-1]_N\big)
\eeq
which is implied by the fact that $\tau+E[-1]$ is a Manin matrix;
see \cite{cfr:ap} for an extensive review on Manin matrices.
Indeed, by \eqref{am} for the image of \eqref{idtrd} under $\vs$ we get
\ben
\tr\ts A^{(N)} \big(\tau-E[-1]_1\big)\dots \big(\tau-E[-1]_N\big)
=\tr\ts A^{(N)} \big(\tau-E^{\tss t}[-1]_1\big)\dots \big(\tau-E^{\tss t}[-1]_N\big),
\een
where we have applied the transposition $t_1\dots t_N$ with respect to all copies of $\End\CC^{N}$
and used the invariance of $A^{(N)}$ under this transposition. Since
$\tau-E^{\tss t}[-1]$ is also a Manin matrix, the resulting
expression coincides with $\cdet\ts \big(\tau-E^{\tss t}[-1]\big)$.
Finally, \eqref{tr} follows from the Newton-type formula connecting the coefficients
of the polynomial in \eqref{deftracepa} with those of \eqref{coldetal}; see
\cite[(3.5)]{cm:ho}.
\epf

With the parameters chosen as in Sec.~\ref{sec:ga}, suppose that $\chi$ vanishes
on the subspace $\n_-\oplus\n_+$ of $\gl_N$ so that we can regard $\chi$
as an element of $\h^*$. Set
\ben
E_{ij}(u)=\sum_{a=1}^\ell \frac{(E_{ij})_a}{u-z_a}-\chi(E_{ij})\in \U(\gl_N)^{\ot\tss\ell}.
\een
Consider the row-determinant $\rdet\big(\di_u+E(u)\big)$ of the matrix
$\di_u+E(u)=\big[\de_{ij}\tss \di_u+E_{ij}(u)\big]$ as a differential operator in $\di_u$
with coefficients in $\U(\gl_N)^{\ot\tss\ell}$.
Furthermore, in accordance with \eqref{hw},
set
\ben
\Ec_{ii}(u)=\sum_{a=1}^\ell \frac{\la_a(E_{ii})}
{u-z_a}-\sum_{j=1}^m\frac{\al_{i_j}(E_{ii})}{u-w_j}-\chi(E_{ii}).
\een

In all the following eigenvalue formulas for the Gaudin Hamiltonians
we will assume that the Bethe ansatz equations \eqref{bae} hold.

\bth\label{thm:rdet}
The eigenvalue of the operator
$\rdet\big(\di_u+E(u)\big)$ on the Bethe vector \eqref{bv} is found by
\ben
\rdet\big(\di_u+E(u)\big)\ts \phi(w_1^{i_1},\dots,w_m^{i_m})
=\big(\di_u+\Ec_{NN}(u)\big)\dots \big(\di_u+\Ec_{11}(u)\big)\ts \phi(w_1^{i_1},\dots,w_m^{i_m}).
\een
\eth

\bpf
To apply Theorem~\ref{thm:eigen},
we will find the image of the polynomial $\cdet\big(\tau+E[-1]\big)$ under the
anti-homomorphism $\Phi$.
We regard $\Phi$ as the map
\ben
\Phi:\U\big(t^{-1}\gl_N[t^{-1}]\big)\ot\CC[\tau]\to \U(\gl_N)^{\ot\tss\ell}\ot \CC[\di_u]
\een
such that $\tau\mapsto\di_u$. Note that by definition of the homomorphism \eqref{psiu}
we have
\ben
\Psi:E[-1]\mapsto -E(u)
\een
and so, by \eqref{cdetim},
\beql{colre}
\Phi:\cdet\ts \big(\tau+E[-1]\big)\mapsto
\cdet\big(\di_u+E^t(u)\big)=\rdet\big(\di_u+E(u)\big).
\eeq

The images of the elements $\vp_i$ under the isomorphism \eqref{hchiaff}
for $\g=\gl_N$ are easy to obtain from \eqref{coldetal}, they are found by
\beql{hchcdet}
\f:\cdet\big(\tau+E[-1]\big)\mapsto \big(\tau+E_{NN}[-1]\big)\dots \big(\tau+E_{11}[-1]\big).
\eeq
Therefore,
\ben
\varrho\circ\f:\cdet\big(\tau+E[-1]\big)
\mapsto
\big(\di_u+\Ec_{NN}(u)\big)\dots \big(\di_u+\Ec_{11}(u)\big)
\een
completing the proof.
\epf

Formula \eqref{hchcdet} can be generalized to get the Harish-Chandra images
of the polynomials \eqref{deftra} and \eqref{deftrh}. We get
\begin{align}
\f:\tr\ts A^{(m)}\big(\tau+E[-1]_1\big)\dots&\big(\tau+E[-1]_m\big)
\mapsto e_m\big(\tau+E_{11}[-1],\dots,\tau+E_{NN}[-1]\big),
\non\\[0.7em]
\f:\tr\ts H^{(m)}\big(\tau+E[-1]_1\big)\dots&\big(\tau+E[-1]_m\big)
\mapsto h_m\big(\tau+E_{11}[-1],\dots,\tau+E_{NN}[-1]\big),
\non
\end{align}
where we use standard noncommutative versions of the complete
and elementary symmetric functions in the ordered variables $x_1,\dots,x_p$
defined by the respective formulas
\begin{align}\label{comp}
h_m(x_1,\dots,x_p)&=\sum_{i_1\leqslant\dots\leqslant i_m}
x_{i_1}\dots x_{i_m},\\
\label{elem}
e_m(x_1,\dots,x_p)&=\sum_{i_1>\dots> i_m}
x_{i_1}\dots x_{i_m}.
\end{align}
The following corollaries can be derived from Theorem~\ref{thm:rdet}
or proved in a similar way with the use of Lemma~\ref{lem:sgn}.

\bco\label{cor:trah}
The eigenvalues of the operators
\ben
\tr\ts A^{(m)}\big(\di_u+E(u)_1\big)\dots\big(\di_u+E(u)_m\big)
\fand
\tr\ts H^{(m)}\big(\di_u+E(u)_1\big)\dots\big(\di_u+E(u)_m\big)
\een
on the Bethe vector \eqref{bv} are found by
respective formulas
\ben
e_m\big(\di_u+\Ec_{11}(u),\dots,\di_u+\Ec_{NN}(u)\big)
\fand
h_m\big(\di_u+\Ec_{11}(u),\dots,\di_u+\Ec_{NN}(u)\big).
\een
\eco

By \cite[Corollary~6.4]{cm:ho} we have
\begin{multline}
\f:\sum_{k=0}^{\infty}z^k\ts \tr\big(\tau+E[-1]\big)^k
\mapsto
\sum_{i=1}^N\Big(1-z\tss\big(\tau+E_{11}[-1]\big)\Big)^{-1}\cdots
\Big(1-z\tss\big(\tau+E_{i\tss i}[-1]\big)\Big)^{-1}\\
{}\times{}\Big(1-z\tss\big(\tau+E_{i-1\ts i-1}[-1]\big)\Big)\cdots
\Big(1-z\tss\big(\tau+E_{11}[-1]\big)\Big),
\non
\end{multline}
where $z$ is an independent variable. So we get the following.

\bco\label{cor:trap}
The eigenvalue of the series
\ben
\sum_{k=0}^{\infty}z^k\ts \tr\big(\di_u+E^t(u)\big)^k
\een
on the Bethe vector \eqref{bv} is found by
the formula
\begin{multline}
\sum_{i=1}^N\Big(1-z\tss\big(\di_u+\Ec_{11}(u)\big)\Big)^{-1}\cdots
\Big(1-z\tss\big(\di_u+\Ec_{i\tss i}(u)\big)\Big)^{-1}\\
{}\times{}\Big(1-z\tss\big(\di_u+\Ec_{i-1\ts i-1}(u)\big)\Big)\cdots
\Big(1-z\tss\big(\di_u+\Ec_{11}(u)\big)\Big).
\non
\end{multline}
\eco

\subsection{Types $B$ and $D$}
\label{subssec:typeb}

Now turn to the orthogonal Lie algebras and
let $\g=\oa_N$ with $N=2n$ or $N=2n+1$. These are simple Lie
algebras of types $D_n$ and $B_n$, respectively.
We use the involution
on the set $\{1,\dots,N\}$ defined by $i'=N-i+1$.
The Lie subalgebra
of $\gl_N$ spanned by the elements $F_{ij}=E_{ij}-E_{j'i'}$
with $i,j\in\{1,\dots,N\}$
is isomorphic to
the orthogonal Lie algebra $\oa_N$.

Denote by $\h$ the Cartan subalgebra of $\oa_N$ spanned by the basis
elements $F_{11},\dots,F_{nn}$.
We have the triangular decomposition $\oa_N=\n_-\oplus\h\oplus\n_+$,
where $\n_-$ and $\n_+$ denote the subalgebras of $\oa_N$ spanned by the
elements
$F_{ij}$ with $i>j$ and by
the elements
$F_{ij}$ with $i<j$, respectively.

We will use the elements
$F_{ij}[r]=F_{ij}\tss t^r$ of the loop algebra $\oa_N[t,t^{-1}]$.
Introduce the elements $F[r]_a$ of the algebra \eqref{tenprka}
by
\beql{matnot}
F[r]_a=\sum_{i,j=1}^{N}
1^{\ot(a-1)}\ot e_{ij}\ot 1^{\ot(m-a)}\ot F_{ij}[r],
\eeq
where $\U$ in \eqref{tenprka} now stands
for the universal enveloping algebra of
$\wh\oa_N\oplus\CC\tau$.

For $1\leqslant a<b\leqslant m$ consider
the operators $P_{a\tss b}$ defined by \eqref{pdef} and introduce
the operators
\ben
Q_{a\tss b}=\sum_{i,j=1}^N 1^{\ot(a-1)}\ot e_{ij}
\ot 1^{\ot(b-a-1)}\ot e_{i'j'}\ot 1^{\ot(m-b)}.
\een
Set
\beql{symo}
S^{(m)}=\frac{1}{m!}
\prod_{1\leqslant a<b\leqslant m}
\Big(1+\frac{P_{a\tss b}}{b-a}-\frac{Q_{a\tss b}}
{N/2+b-a-1}\Big),
\eeq
where the product is taken in the lexicographic order
on the pairs $(a,b)$.
The element \eqref{symo} is the image of the
symmetrizer in the Brauer algebra $\Bc_m(N)$ under its
action on the vector space $(\CC^N)^{\ot m}$. In particular,
for any $1\leqslant a<b\leqslant m$ for the operator $S^{(m)}$ we have
\beql{sqpo}
S^{(m)}\ts Q_{a\tss b}=Q_{a\tss b}\ts S^{(m)}=0
\Fand S^{(m)}\ts P_{a\tss b}=P_{a\tss b}\ts S^{(m)}=S^{(m)}.
\eeq
The symmetrizer admits a few
other equivalent expressions which are reproduced in \cite{m:ff}.

We will use the notation
\beql{galm}
\ga_m(\om)=\frac{\om+m-2}{\om+2\tss m-2}
\eeq
and define the elements $\vp^{}_{m\tss a}\in\U\big(t^{-1}\oa_N[t^{-1}]\big)$
by the expansion
\beql{deftr}
\ga_m(N)\ts\tr\ts S^{(m)} \big(\tau+F[-1]_1\big)\dots \big(\tau+F[-1]_m\big)
=\vp^{}_{m\tss0}\ts\tau^m+\vp^{}_{m\tss1}\ts\tau^{m-1}
+\dots+\vp^{}_{m\tss m},
\eeq
where the trace is taken
over all $m$ copies of $\End\CC^{N}$.
By the main result of \cite{m:ff}, all coefficients $\vp^{}_{m\tss a}$
belong to the Feigin--Frenkel center $\z(\wh\oa_N)$.
In the even orthogonal case $\g=\oa_{2n}$ there is an additional element
$\vp^{\tss\prime}_n=\Pf\ts\wt F[-1]$ of the center defined as
the (noncommutative) Pfaffian of the skew-symmetric matrix $\wt F[-1]=\big[\wt F_{ij}[-1]\big]$,
\beql{genpf}
\Pf\ts\wt F[-1]=\frac{1}{2^nn!}\sum_{\si\in\Sym_{2n}}\sgn\si\cdot
\wt F_{\si(1)\ts\si(2)}[-1]\dots \wt F_{\si(2n-1)\ts\si(2n)}[-1],
\eeq
where $\wt F_{ij}[-1]=F_{ij'}[-1]$. The family
$\vp^{}_{2\tss 2},\vp^{}_{4\tss 4},\dots,
\vp^{}_{2n\ts 2n}$
is a complete set of Segal--Sugawara vectors
for $\oa_{2n+1}$, whereas
$\vp^{}_{2\tss 2},\vp^{}_{4\tss 4},\dots,\vp^{}_{2n-2\ts 2n-2},
\vp^{\tss\prime}_n$
is a complete set of Segal--Sugawara vectors
for $\oa_{2n}$.

We extend the involution \eqref{sgn} to the algebra
$\U\big(t^{-1}\oa_N[t^{-1}]\big)\ot\CC[\tau]$ with the action on $\CC[\tau]$
as the identity map.

\ble\label{lem:sgnbcd}
The element \eqref{deftr} is stable under $\vs$. Moreover,
in type $D_n$ we have
\beql{pfaims}
\vs:\Pf\ts\wt F[-1]\mapsto (-1)^n\ts\Pf\ts\wt F[-1].
\eeq
\ele

\bpf
The same argument as in the proof of Lemma~\ref{lem:sgn} shows that
the image of \eqref{deftr} under the involution $\vs$ equals
\beql{deftrbd}
\ga_m(N)\ts\tr\ts S^{(m)} \big(\tau-F[-1]_1\big)\dots \big(\tau-F[-1]_m\big).
\eeq
Indeed, this is implied by \eqref{sqpo} and the commutation relations
\ben
F[r]_a\ts F[s]_b-F[s]_b\ts F[r]_a
=(P_{a\tss b}-Q_{a\tss b})\ts F[r+s]_b-F[r+s]_b\ts (P_{a\tss b}-Q_{a\tss b})
\een
for $a<b$. By applying
the simultaneous transpositions
$e_{ij}\mapsto e_{j'i'}$ to all $m$ copies of $\End\CC^N$ we conclude
that \eqref{deftrbd} coincides with \eqref{deftr} because
this transformation takes
each factor $\tau-F[-1]_a$ to $\tau+F[-1]_a$ whereas
the operator $S^{(m)}$ stays invariant. Relation \eqref{pfaims}
is immediate from \eqref{genpf}.
\epf

By the main results of \cite{mm:yc}, the image of the polynomial \eqref{deftr} under
the isomorphism \eqref{hchiaffe}
is given by the formula{\tss\rm:}
\ben
h_m\big(\tau+F_{11}[-1],\dots,
\tau+F_{n\tss n}[-1],\tau-F_{n\tss n}[-1],\dots
\tau-F_{11}[-1]\big),
\een
for type $B_n$ and by
\begin{multline}
{\textstyle \frac{1}{2}}\ts h_m\big(\tau+F_{11}[-1],\dots,
\tau+F_{n-1\ts n-1}[-1],\tau-F_{n\tss n}[-1],\dots
\tau-F_{11}[-1]\big)\\[0.8em]
{}{\qquad\qquad\qquad\qquad}\quad+{\textstyle \frac{1}{2}}\ts h_m\big(\tau+F_{11}[-1],\dots,
\tau+F_{n\tss n}[-1],\tau-F_{n-1\ts n-1}[-1],\dots
\tau-F_{11}[-1]\big),
\non
\end{multline}
for type $D_n$. The latter sum can also be written in the form
\begin{multline}
h_m\big(\tau+F_{11}[-1],\dots,
\tau+F_{n\tss n}[-1],\tau-F_{n\tss n}[-1],\dots
\tau-F_{11}[-1]\big)\\[0.8em]
{}-\sum_{k+l=m-1}h_k\big(\tau+F_{11}[-1],\dots,
\tau+F_{n\tss n}[-1]\big)\ts\tau\ts h_{\tss l}\big(\tau-F_{n\tss n}[-1],\dots,
\tau-F_{11}[-1]\big).
\non
\end{multline}
Furthermore, the image
of the element $\vp^{\tss\prime}_n$ in type $D_n$ is given by
\beql{pfaim}
\big(F_{11}[-1]-\tau\big)\dots \big(F_{n\tss n}[-1]-\tau\big)\ts 1,
\eeq
where $\tau$ is understood as the differentiation operator
so that $\tau\ts 1=0$; see also \cite{r:nf} for a direct calculation
of the Harish-Chandra image of $\vp^{\tss\prime}_n$.

Choose parameters as in Sec.~\ref{sec:ga} and suppose that $\chi$ vanishes
on the subspace $\n_-\oplus\n_+$ of $\oa_N$ so that we can regard $\chi$
as an element of $\h^*$. Set
\ben
F_{ij}(u)=\sum_{a=1}^\ell \frac{(F_{ij})_a}{u-z_a}-\chi(F_{ij})\in \U(\oa_N)^{\ot\tss\ell}.
\een
In accordance with \eqref{hw}
set
\ben
\Fc_{ii}(u)=\sum_{a=1}^\ell \frac{\la_a(F_{ii})}
{u-z_a}-\sum_{j=1}^m\frac{\al_{i_j}(F_{ii})}{u-w_j}-\chi(F_{ii}).
\een
In the case $\g=\oa_{2n}$ define the operator
\beql{genpfu}
\Pf\ts\wt F(u)=\frac{1}{2^nn!}\sum_{\si\in\Sym_{2n}}\sgn\si\cdot
\wt F_{\si(1)\ts\si(2)}(u)\dots \wt F_{\si(2n-1)\ts\si(2n)}(u),
\eeq
where $\wt F_{ij}(u)=F_{ij'}(u)$.
As before, we will assume that the Bethe ansatz equations \eqref{bae} hold.

\bth\label{thm:symbd}
The eigenvalue of the operator
\beql{deftrim}
\ga_m(N)\ts\tr\ts S^{(m)} \big(\di_u+F(u)_1\big)\dots \big(\di_u+F(u)_m\big)
\eeq
on the Bethe vector \eqref{bv} is found by
\ben
h_m\big(\di_u+\Fc_{11}(u),\dots,
\di_u+\Fc_{n\tss n}(u),\di_u-\Fc_{n\tss n}(u),\dots
\di_u-\Fc_{11}(u)\big)
\een
for type $B_n$, and by
\begin{multline}
{\textstyle \frac{1}{2}}\ts h_m\big(\di_u+\Fc_{11}(u),\dots,
\di_u+\Fc_{n-1\ts n-1}(u),\di_u-\Fc_{n\tss n}(u),\dots
\di_u-\Fc_{11}(u)\big)\\[0.8em]
{}{\quad\quad\qquad\qquad}\quad+{\textstyle \frac{1}{2}}\ts h_m\big(\di_u+\Fc_{11}(u),\dots,
\di_u+\Fc_{n\tss n}(u),\di_u-\Fc_{n-1\ts n-1}(u),\dots
\di_u-\Fc_{11}(u)\big)
\non
\end{multline}
for type $D_n$.
Moreover, the eigenvalue of the
operator $\Pf\ts\wt F(u)$ in type $D_n$ is given by
\beql{pfaimu}
\big(\Fc_{11}(u)-\di_u\big)\dots \big(\Fc_{n\tss n}(u)-\di_u\big)\ts 1.
\eeq
\eth

\bpf
We apply Theorem~\ref{thm:eigen} again
and regard $\Phi$ as the map
\ben
\Phi:\U\big(t^{-1}\oa_N[t^{-1}]\big)\ot\CC[\tau]\to \U(\oa_N)^{\ot\tss\ell}\ot \CC[\di_u]
\een
such that $\tau\mapsto\di_u$. By the definition of the homomorphism \eqref{psiu}
we have
\ben
\Psi:F[-1]\mapsto -F(u).
\een
Hence, using the equivalent formula \eqref{deftrbd} for the polynomial \eqref{deftr}
we find that its image under $\Phi$ coincides with
the operator \eqref{deftrim}.
The proof of the first part of the theorem is completed by using
the formulas for the images of \eqref{deftr} under
the respective isomorphisms \eqref{hchiaffe} recalled above.
Finally, by Lemma~\ref{lem:sgnbcd}, in type $D_n$,
\ben
\Phi:\Pf\ts\wt F[-1]\mapsto \Pf\ts\wt F(u)
\een
so that the last claim follows by using formula \eqref{pfaim} for the
image of $\Pf\ts\wt F[-1]$ under the isomorphism \eqref{hchiaff}.
\epf

\bco\label{cor:genf}
The eigenvalue of the generating function
\beql{deftrimgf}
\Bigg(\sum_{m=0}^\infty (-z)^m\ts\ga_m(N)\ts\tr\ts S^{(m)}
\big(\di_u+F(u)_1\big)\dots \big(\di_u+F(u)_m\big)\Bigg)^{-1}
\eeq
on the Bethe vector \eqref{bv} is found by
\ben
\Big(1+\big(\di_u-\Fc_{11}(u)\big)z\Big)\dots\Big(1+\big(\di_u-\Fc_{n\tss n}(u)\big)z\Big)
\Big(1+\big(\di_u+\Fc_{n\tss n}(u)\big)z\Big)\dots \Big(1+\big(\di_u+\Fc_{11}(u)\big)z\Big)
\een
for type $B_n$ and by
\begin{multline}
\Big(1+\big(\di_u-\Fc_{11}(u)\big)z\Big)\dots\Big(1+\big(\di_u-\Fc_{n\tss n}(u)\big)z\Big)
\Big(1+\di_u\tss z\Big)^{-1}\\[0.8em]
{}\times\Big(1+\big(\di_u+\Fc_{n\tss n}(u)\big)z\Big)\dots \Big(1+\big(\di_u+\Fc_{11}(u)\big)z\Big)
\non
\end{multline}
for type $D_n$.
\eco

\subsection{Type $C$}
\label{subssec:typec}

We identify the symplectic Lie algebra $\g=\spa_{2n}$
with the Lie subalgebra of $\gl_{2n}$ spanned by the elements
$F_{ij}=E_{ij}-\ve_i\tss\ve_j\tss E_{j'i'}$ with
$i,j\in\{1,\dots,2n\}$, where $i'=2n-i+1$ and
$\ve_i=1$ for $i=1,\dots,n$ and
$\ve_i=-1$ for $i=n+1,\dots,2n$.

Denote by $\h$ the Cartan subalgebra of $\spa_{2n}$ spanned by the basis
elements $F_{11},\dots,F_{nn}$.
We have the triangular decomposition $\spa_{2n}=\n_-\oplus\h\oplus\n_+$,
where $\n_-$ and $\n_+$ denote the subalgebras of $\spa_{2n}$ spanned by the
elements
$F_{ij}$ with $i>j$ and by
the elements
$F_{ij}$ with $i<j$, respectively.

We will use the elements
$F_{ij}[r]=F_{ij}\tss t^r$ of the loop algebra $\spa_{2n}[t,t^{-1}]$.
Introduce the elements $F[r]_a$ of the algebra \eqref{tenprka}
by
\beql{matnotc}
F[r]_a=\sum_{i,j=1}^{2n}
1^{\ot(a-1)}\ot e_{ij}\ot 1^{\ot(m-a)}\ot F_{ij}[r],
\eeq
where $\U$ in \eqref{tenprka} now stands
for the universal enveloping algebra of
$\wh\spa_{2n}\oplus\CC\tau$.

For $1\leqslant a<b\leqslant m$ consider
the operators $P_{a\tss b}$ defined by \eqref{pdef} and introduce
the operators
\ben
Q_{a\tss b}=\sum_{i,j=1}^{2n} \ve_i\tss\ve_j\ts 1^{\ot(a-1)}\ot e_{ij}
\ot 1^{\ot(b-a-1)}\ot e_{i'j'}\ot 1^{\ot(m-b)}.
\een
For $1\leqslant m\leqslant n$ set
\beql{symsp}
S^{(m)}=\frac{1}{m!}
\prod_{1\leqslant a<b\leqslant m}
\Big(1-\frac{P_{a\tss b}}{b-a}-\frac{Q_{a\tss b}}
{n-b+a+1}\Big),
\eeq
where the product is taken in the lexicographic order
on the pairs $(a,b)$.
The element \eqref{symsp} is the image of the
symmetrizer in the Brauer algebra $\Bc_m(-2n)$ under its
action on the vector space $(\CC^{2n})^{\ot m}$.
Use the notation \eqref{galm} to introduce the polynomial in $\tau$ by
\beql{deftrc}
\ga_m(-2n)\ts\tr\ts S^{(m)} \big(\tau+F[-1]_1\big)\dots \big(\tau+F[-1]_m\big)
=\vp^{}_{m\tss0}\ts\tau^m+\vp^{}_{m\tss1}\ts\tau^{m-1}
+\dots+\vp^{}_{m\tss m},
\eeq
where the trace is taken
over all $m$ copies of $\End\CC^{2n}$.
By the results of \cite{m:ff}, the values of $m$ in
\eqref{deftrc} can be extended to the range
$1\leqslant m\leqslant 2n$ (and, in fact, for $m=2n+1$ as well) to get
a well-defined polynomial in $\tau$.
Moreover, the family
$\vp^{}_{2\tss 2},\vp^{}_{4\tss 4},\dots,
\vp^{}_{2n\ts 2n}$
is a complete set of Segal--Sugawara vectors
for $\spa_{2n}$.

Extend the involution \eqref{sgn} to the algebra
$\U\big(t^{-1}\spa_{2n}[t^{-1}]\big)\ot\CC[\tau]$ with the action on $\CC[\tau]$
as the identity map.

\ble\label{lem:sgnc}
The element \eqref{deftrc} is stable under $\vs$.
\ele

\bpf
The proof is the same as for Lemma~\ref{lem:sgnbcd}, which also provides
an equivalent formula
\beql{deftreq}
\ga_m(-2n)\ts\tr\ts S^{(m)} \big(\tau-F[-1]_1\big)\dots \big(\tau-F[-1]_m\big)
\eeq
for the polynomial \eqref{deftrc}.
\epf

By the main result of \cite{mm:yc}, the image of the polynomial \eqref{deftrc}
with $1\leqslant m\leqslant 2n+1$ under
the isomorphism \eqref{hchiaffe}
is given by the formula{\tss\rm:}
\ben
e_m\big(\tau+F_{11}[-1],\dots,
\tau+F_{n\tss n}[-1],\tau,\tau-F_{n\tss n}[-1],\dots
\tau-F_{11}[-1]\big),
\een
where we use notation \eqref{elem}.

With parameters chosen as in Sec.~\ref{sec:ga},
suppose that $\chi$ vanishes
on the subspace $\n_-\oplus\n_+$ of $\spa_{2n}$ so that we can regard $\chi$
as an element of $\h^*$. Set
\ben
F_{ij}(u)=\sum_{a=1}^\ell \frac{(F_{ij})_a}{u-z_a}-\chi(F_{ij})\in \U(\spa_{2n})^{\ot\tss\ell}.
\een
In accordance with \eqref{hw}
set
\ben
\Fc_{ii}(u)=\sum_{a=1}^\ell \frac{\la_a(F_{ii})}
{u-z_a}-\sum_{j=1}^m\frac{\al_{i_j}(F_{ii})}{u-w_j}-\chi(F_{ii}).
\een
As before, we will assume that the Bethe ansatz equations \eqref{bae} hold.

\bth\label{thm:symc}
For any $1\leqslant m\leqslant 2n+1$ the eigenvalue of the operator
\beql{deftrimc}
\ga_m(-2n)\ts\tr\ts S^{(m)} \big(\di_u+F(u)_1\big)\dots \big(\di_u+F(u)_m\big)
\eeq
on the Bethe vector \eqref{bv} is found by
\ben
e_m\big(\di_u+\Fc_{11}(u),\dots,
\di_u+\Fc_{n\tss n}(u),\di_u,\di_u-\Fc_{n\tss n}(u),\dots
\di_u-\Fc_{11}(u)\big).
\een
\eth

\bpf
This is derived from Theorem~\ref{thm:eigen} and Lemma~\ref{lem:sgnc}
as in the proof of Theorem~\ref{thm:symbd}.
\epf

\subsection{Connection with the results of \cite{mtv:be} and \cite{mv:cp}}
\label{subssec:comp}

Theorem~\ref{thm:rdet} was previously proved
in \cite{mtv:be} is a slightly different form; see Theorem~9.2 there.
We will make a connection between these results by showing that one is obtained
from the other by using an automorphism of the current algebra.
The notation of \cite{mtv:be} corresponds to ours (we used the settings of \cite{ffr:gm}
and \cite{fft:gm}) as follows.
The highest weights
$\La_k=(\La^1_k,\dots,\La^N_k)$ correspond to our $\la_k$ so that $\La^i_k=\la_k(E_{ii})$;
the evaluation parameters $z_i$ are the same.
The diagonal matrix $K=\diag\ts[K_1,\dots,K_N]$ corresponds to our element $-\chi$ so that
$K_i=-\chi(E_{ii})$. Finally, the collection of nonnegative integers
$\xi=(\xi^1,\dots,\xi^{N-1})$ gives rise to
our multiset of simple roots $\al_{i_j}$ where $\al_l=\ve_l-\ve_{l+1}$ occurs $\xi^l$ times for
each $l=1,\dots,N-1$.
The corresponding variables
$t^1_1,\dots,t^1_{\xi^1},\dots,t^{N-1}_1,\dots,t^{N-1}_{\xi^{N-1}}$
are then respectively identified with our parameters
$w_1,\dots,w_m$ with $m=|\xi|$. The coroots $\check\al_l$ coincide with the elements
$E_{l\tss l}-E_{l+1\tss l+1}$
so that the Bethe ansatz equations (9.3) in \cite{mtv:be} turn into
\eqref{bae}.
Using this correspondence between the settings, we can now state \cite[Theorem~9.2]{mtv:be}
in our notation as the relation
\ben
\cdet\big(\di_u-E(u)\big)\ts \phi(w_1^{i_1},\dots,w_m^{i_m})
=\big(\di_u-\Ec_{11}(u)\big)\dots \big(\di_u-\Ec_{NN}(u)\big)\ts \phi(w_1^{i_1},\dots,w_m^{i_m})
\een
for the eigenvalue of the operator $\cdet\big(\di_u-E(u)\big)$ on
the Bethe vector \eqref{bv}.
This relation is implied by Theorem~\ref{thm:rdet} by
twisting the action of $\U(\gl_N)$ on each Verma module $M_{\la_k}$
by the automorphism $E_{ij}\mapsto -E_{j'i'}$, where $i'=N-i+1$.
The automorphism takes $\rdet\big(\di_u+E(u)\big)$ to $\cdet\big(\di_u-E(u)\big)$
and $\Ec_{i\tss i}(u)$ to $-\Ec_{i'i'}(u)$.

We also make a connection of
Theorems~\ref{thm:rdet}, \ref{thm:symbd} and \ref{thm:symc}
with formulas for universal differential operators corresponding to
populations of critical points of the master functions associated with
flag varieties; see \cite{mv:cp}. With the recalled above notation of \cite{mtv:be},
we follow \cite{mv:cp} to
introduce polynomials in type $A$,
\ben
T_a(u)=\prod_{k=1}^{\ell}(u-z_k)^{\La_k^a},\qquad a=1,\dots,N,
\een
and
\ben
y_a(u)=\prod_{p=1}^{\xi^a}(u-t^a_p),\qquad a=1,\dots,N-1.
\een
Then the eigenvalue of the Bethe vector in Theorem~\ref{thm:rdet} with $\chi=0$
coincides with the differential operator
\beql{diffa}
\prod_{a=1,\dots,N}^{\longleftarrow}\Bigg(\di_u+\ln'\tss\frac{T_a(u)\ts y_{a-1}(u)}{y_a(u)}\Bigg)
\eeq
rewritten in our notation,
where we set $y_0(u)=y_N(u)=1$; see \cite[Sec.~5.2]{mv:cp}.

Using a similar notation, in type $B_n$ set
\beql{tau}
T^B_a(u)=\prod_{k=1}^{\ell}(u-z_k)^{\La_k^a}\Fand
y^B_a(u)=\prod_{p=1}^{\xi^a}(u-t^a_p),
\qquad a=1,\dots,n.
\eeq
Then the coefficient of $z^{2n}$
in the eigenvalue in type $B_n$ (see Corollary~\ref{cor:genf})
coincides with \eqref{diffa}, if we take $N=2n$ and set
\ben
y_a(u)=y_{2n-a}(u)=y^B_a(u)\qquad\text{for}\quad a=1,\dots,n
\een
and
\ben
T_a(u)=T_{2n-a+1}(u)^{-1}=T^B_a(u)\qquad\text{for}\quad a=1,\dots,n;
\een
cf. \cite[Sec.~7.1]{mv:cp}. In type $C_n$, introducing $T^C_a(u)$ and $y^C_a(u)$
for $a=1,\dots,n$ as in \eqref{tau}, we find that
the eigenvalue of the operator with $m=2n+1$ in
Theorem~\ref{thm:symc} is given by \eqref{diffa} with $N=2n+1$,
where we set
\ben
y_a(u)=y_{2n-a+1}(u)=\begin{cases} y^C_a(u)\qquad&\text{for}\quad a=1,\dots,n-1\\
                                   y^C_a(u)^2\qquad&\text{for}\quad a=n
                                   \end{cases}
\een
and
\ben
T_a(u)=T_{2n-a+2}(u)^{-1}=\begin{cases} T^C_a(u)\qquad&\text{for}\quad a=1,\dots,n\\
                                   1\qquad&\text{for}\quad a=n+1;
                                   \end{cases}
\een
cf. \cite[Sec.~7.2]{mv:cp}.

\section{From $q$-characters to classical $\Wc$-algebras}
\label{sec:cwa}
\setcounter{equation}{0}

The Harish-Chandra images of the Segal--Sugawara elements
\eqref{deftr} and \eqref{deftrc} in types $B$, $C$ and $D$
were calculated in \cite{mm:yc} by taking a classical limit of certain
Yangian characters (or $q$-characters). Our goal in this section
is to prove general results providing a connection between
the rings of $q$-characters and the corresponding classical
$\Wc$-algebras. We will rely on the original work \cite{fr:qc}
for the basic definitions and properties of the $q$-characters; see also \cite{fm:cq}.
However, we will
use an equivalent additive version of the character ring as in \cite{nn:pt}
and indicate the connection between the notation in Remarks~\ref{rem:notla}
and \ref{rem:notlbcd} below. Although this version
can be introduced independently via the Yangian representation theory,
we will not make a direct use of the Yangians which will only appear
in the notation $\Rep\Y(\g)$ for the ring of characters; cf. \cite{mm:yc}.

The screening operators for classical $\Wc$-algebras
are constructed as limits of certain intertwiners between
$\wh\g_{\ka}$-modules at a level $\ka$, as $\ka\to -h^{\vee}$;
see \cite[Ch.~7]{f:lc}. They can also be obtained by applying
a Chevalley-type theorem to the $\Wc$-algebras defined in the context
of classical Hamiltonian reduction; see, e.g., \cite{mr:cw}.
It was conjectured in \cite{fr:qc} and proved in \cite{fm:cq}, that the ring of characters
can be defined as the intersection of the kernels of the screening operators.
We will apply a classical limit procedure
to derive the screening operators characterizing elements of the $\Wc$-algebra;
cf. \cite[Sec.~8]{fr:qc}. The main result of \cite{mm:yc} will play
an important role in the proof of the
surjectivity of the procedure.

\subsection{Type $A$}
\label{subssec:scra}

Introduce the algebra of polynomials
\ben
\Lc=\CC[\la_i(a)\ts |\ts i=1,\dots,N,\ts a\in\CC]
\een
in the variables $\la_i(a)$. For every $i\in\{1,\dots,N-1\}$
consider the free left $\Lc$-module $\wt\Lc_i$ with the
generators $\si_i(a)$, where $a$ runs over $\CC$ and denote by $\Lc_i$
its quotient by the relations
\beql{rellasi}
\la_i(a)\ts \si_i(a)=\la_{i+1}(a)\ts \si_i(a+1),\qquad a\in\CC.
\eeq
Define the linear operator
$\wt S_i:\Lc\to \wt\Lc_i$
by the formula
\beql{defsi}
\wt S_i:\la_j(a)\mapsto\begin{cases} \la_i(a)\ts\si_i(a)\qquad&\text{for}\quad j=i\\
-\la_{i+1}(a)\ts\si_i(a+1)\qquad&\text{for}\quad j=i+1\\
0\qquad&\text{for}\quad j\ne i,i+1
\end{cases}
\eeq
and the Leibniz rule
\beql{lr}
\wt S_i(AB)=B\tss\wt S_i(A)+A\tss\wt S_i(B).
\eeq
Now the {\em $i$-th screening operator}
\ben
S_i:\Lc\to \Lc_i
\een
is defined as the composition of $\wt S_i$ and the projection $\wt\Lc_i\to \Lc_i$.

In accordance with \cite[Theorem~5.1]{fm:cq}, we can define the subalgebra $\Rep\Y(\gl_N)$
of Yangian characters in $\Lc$ as the intersection of kernels of the screening
operators:
\ben
\Rep\Y(\gl_N)=\bigcap_{i=1}^{N-1} \ts\ker S_i.
\een

\bre\label{rem:notla}
Our variables $\la_i(a)$ and $\si_i(a)$ correspond to $\La_{i,q^{2a}}$ and $S_{i,q^{2a+i-1}}$
from \cite{fr:qc},
respectively; cf. \cite{nn:pt}.
\qed
\ere

Now we recall the definition of the classical $\Wc$-algebra $\Wc(\gl_N)$
via screening operators as in \cite[Sec.~8.1]{f:lc}; see also \cite{mm:yc} and \cite{mr:cw}.
With the notation as in Sec.~\ref{subssec:typea}, we will regard $\U(\wh\h_-)$
as the algebra of polynomials in the variables $E_{ii}[r]$
with $i=1,\dots,N$ and $r<0$. The screening operators
\ben
V_i:\U(\wh\h_-)\to \U(\wh\h_-),\qquad i=1,\dots,N-1
\een
are defined by
\ben
V_i=\sum_{r=0}^{\infty} V_{i\ts [r]}\ts
\Bigg(\frac{\di}{\di E_{ii}[-r-1]}-\frac{\di}{\di E_{i+1\ts i+1}[-r-1]}\Bigg),
\een
where the coefficients $V_{i\ts [r]}$ are found from the expansion
of a formal generating function in a variable $z$,
\ben
\sum_{r=0}^{\infty} V_{i\ts [r]}\ts z^r=\exp\ts\sum_{m=1}^{\infty}
\frac{E_{ii}[-m]-E_{i+1\ts i+1}[-m]}{m}\ts z^m.
\een
The {\em classical $\Wc$-algebra} $\Wc(\gl_N)$ is a subalgebra
of $\U(\wh\h_-)$ defined as the intersection of kernels of the screening
operators:
\ben
\Wc(\gl_N)=\bigcap_{i=1}^{N-1} \ts\ker V_i.
\een

We will now construct a map $\gr:\Rep\Y(\gl_N)\to \Wc(\gl_N)$ and describe its properties.
First, embed $\Lc$ into the algebra of formal power series $\CC[[\la^{(r)}_i]]$
in variables $\la^{(r)}_i$ with $i=1,\dots,N$ and $r=0,1,\dots$ by setting
\beql{embedla}
\la_i(a)\mapsto \sum_{r=0}^{\infty}\frac{\la^{(r)}_i}{r!}\ts a^r.
\eeq
Identify the formal power series in the $\la^{(r)}_i$ with those in new variables
$\mu^{(r)}_i$ defined by
\beql{lamu}
\la_i^{(0)}=1+\mu_i^{(0)}\Fand \la_i^{(r)}=\mu_i^{(r)}\quad\text{for}\quad r\geqslant 1.
\eeq
Define the degrees of the new variables by $\deg\mu_i^{(r)}=-r-1$.
Given $A\in\Lc$, consider the corresponding element $\CC[[\mu^{(r)}_i]]$
and take its homogeneous component $\overline A$ of the maximum degree.
This component is a polynomial in the variables $\mu^{(r)}_i$
and so we have a map
\beql{gr}
\gr: \Lc\to \CC[\mu^{(r)}_i],\qquad A\mapsto\overline A.
\eeq
Note its property which is immediate from the definition:
\beql{grab}
\gr(AB)=\gr(A)\ts\gr(B).
\eeq

In the following proposition we identify
$\U(\wh\h_-)$ with the algebra of polynomials $\CC[\mu^{(r)}_i]$ via the isomorphism
$E_{ii}[-r-1]\mapsto \mu^{(r)}_i/r!$.

\bpr\label{prop:grtypea}
The image of the restriction of the map \eqref{gr} to the subalgebra
of characters $\Rep\Y(\gl_N)$ is contained in $\Wc(\gl_N)$ and so it defines a map
\ben
\gr:\Rep\Y(\gl_N)\to \Wc(\gl_N).
\een
Moreover, any homogeneous element
of $\Wc(\gl_N)$ is contained in the image of $\gr$.
\epr

\bpf
Similar to \eqref{embedla},
introduce variables $\si^{(r)}_i$ by the expansion
\beql{embedsi}
\si_i(a)\mapsto \sum_{r=0}^{\infty}\frac{\si^{(r)}_i}{r!}\ts a^r
\eeq
and set $\deg\si^{(r)}_i=-r-1$.
Regarding $a$ as a formal variable in \eqref{rellasi} and \eqref{defsi},
write the screening operators in terms of the variables $\mu^{(r)}_i$.
Explicitly, for $i=1,\dots,N-1$ define
operators
\beql{scirc}
S^{\tss\circ}_i:\CC[[\mu^{(r)}_j]]\to \CC[[\mu^{(r)}_j,\si^{(r)}_i]]\tss/\hspace{-3pt}\thicksim,
\eeq
where the target space is the quotient of $\CC[[\mu^{(r)}_j,\si^{(r)}_i]]$
by the relations \eqref{rellasi} written in terms of the $\mu^{(r)}_i$ with $a$ understood
as a variable. Set
\ben
S^{\tss\circ}_i:\mu^{(0)}_j\mapsto\begin{cases} \big(1+\mu^{(0)}_i\big)\ts \si^{(0)}_i
\qquad&\text{for}\quad j=i\\[0.2em]
-\big(1+\mu^{(0)}_{i+1}\big)\ts {\displaystyle \sum_{k\geqslant 0}}\dfrac{{}\ts\si^{(k)}_i}{k!}
\qquad&\text{for}\quad j=i+1\\
0\qquad&\text{for}\quad j\ne i,i+1
\end{cases}
\een
and
\ben
S^{\tss\circ}_i:\mu^{(r)}_j\mapsto \di^{\ts r}\ts\Big(S^{\tss\circ}_i\big(\mu^{(0)}_j\big)\Big),
\qquad r\geqslant 1,
\een
where the derivation $\di$ acts on the variables by the rule
\ben
\di:\mu^{(r)}_j\mapsto \mu^{(r+1)}_j,\qquad \si^{(r)}_j\mapsto \si^{(r+1)}_j,\qquad r\geqslant 0.
\een
The action of $S^{\tss\circ}_i$ then extends to the entire algebra $\CC[[\mu^{(r)}_j]]$
via the Leibniz rule as in \eqref{lr}.

Now suppose that $A\in \Rep\Y(\gl_N)$ so that $S_i\tss A=0$ for all $i=1,\dots,N-1$.
Denote by $A^{\circ}$ the corresponding
element of $\CC[[\mu^{(r)}_j]]$. By the definition of the operators
$S^{\tss\circ}_i$, their restriction to the subalgebra $\Lc$ coincides
with the action of the respective operators $S_i$. Therefore,
$S^{\tss\circ}_i\tss A^{\circ}=0$.
Taking the top degree component $\overline A$ of $A^{\circ}$
we can write
\ben
S^{\tss\circ}_i\tss A^{\circ}=\overline S_i\tss \overline A+\text{lower degree terms},
\een
where the operator $\overline S_i$ is
given by
\beql{barsi}
\overline S_i:\mu^{(r)}_j\mapsto\begin{cases} \si^{(r)}_i\qquad&\text{for}\quad j=i\\
-\si^{(r)}_{i}\qquad&\text{for}\quad j=i+1\\
0\qquad&\text{for}\quad j\ne i,i+1.
\end{cases}
\eeq
On the other hand, relations \eqref{rellasi} give
\beql{mua}
\mu_i(a)\ts \si_i(a)=\big(1+\mu_{i+1}(a)\big)\ts
\sum_{k=0}^{\infty}\frac{\si^{(k)}_i(a)}{k!}-\si_i(a),
\eeq
where $\si^{(k)}_i(a)$ is defined as the $k$-th derivative over $a$ from \eqref{embedsi}
and
\ben
\mu_j(a)=\sum_{r=0}^{\infty}\frac{\mu^{(r)}_j}{r!}\ts a^r.
\een
Regarding $a$ as a variable, we get from \eqref{mua} a sequence of relations
by comparing the coefficients of the same powers of $a$.
The top degree components in these relations
are homogeneous relations which can be written in terms of generating functions
in the form
\ben
\si'_i(z)=\big(\mu_i(z)-\mu_{i+1}(z)\big)\ts \si_i(z)
\een
so that for the images under
$\overline S_i$ we have
\ben
\overline S_i:\mu_i(z)\mapsto \exp\int\big(\mu_i(z)-\mu_{i+1}(z)\big)\ts dz,\qquad
\mu_{i+1}(z)\mapsto -\exp\int\big(\mu_i(z)-\mu_{i+1}(z)\big)\ts dz,
\een
and $\overline S_i:\mu_j(z)\mapsto 0$ for $j\ne i,i+1$. However, this coincides
with the action of the operator $V_i$ on the series
\ben
\mu_k(z)=\sum_{r=0}^{\infty} E_{kk}[-r-1]\ts z^r,\qquad k=1,\dots,N.
\een
Thus, we may conclude that if an element $A\in\Lc$ is annihilated
by all operators $S_i$, then its image $\overline A$ under the map \eqref{gr} is
annihilated by all operators $V_i$ completing the proof of the first part
of the proposition.

The second part follows from \cite{mm:yc}, where generators of the algebra
$\Wc(\gl_N)$ were obtained as images of certain elements of $\Lc$
under the map $\gr$.
\epf

\subsection{Types $B$, $C$ and $D$}
\label{subssec:scrbcd}

We let $\g$ denote the orthogonal Lie algebra $\oa_N$ (with $N=2n$ or $N=2n+1$)
or the symplectic Lie algebra $\spa_N$ (with $N=2n$).
Introduce a parameter $\ka$ by
$\kappa=N/2-1$ in the orthogonal case and
$\kappa=N/2+1$ in the symplectic case. As before, we set $i'=N-i+1$.

Consider the algebra of polynomials in variables $\la_i(a)$ with
$i=1,\dots,N$ and $a\in\CC$ and denote by $\Lc=\Lc(\g)$ its quotient by
the relations
\beql{larel}
\la_i(a+\ka-i)\ts\la_{i'}(a)=\la_{i+1}(a+\ka-i)\ts\la_{(i+1)'}(a),\qquad a\in\CC,
\eeq
for $i=0,1,\dots,n-1$ if $\g=\oa_{2n}$ or $\spa_{2n}$,
and for $i=0,1,\dots,n$ if $\g=\oa_{2n+1}$, where $\la_0(a)=\la_{0'}(a)=1$.

For $i=1,\dots,n$
consider the free left $\Lc$-module $\wt\Lc_i$ with the
generators $\si_i(a)$, where $a$ runs over $\CC$ and denote by $\Lc_i$
its quotient by the relations
\beql{rellasibcd}
\la_i(a)\ts \si_i(a)=\la_{i+1}(a)\ts \si_i(a+1),\qquad i=1,\dots,n-1,\quad a\in\CC,
\eeq
together with
\beql{quo}
\begin{alignedat}{2}
\la_n(a)\ts \si_n(a)&=\la_{n+1}(a)\ts \si_n(a+1/2),\quad\qquad&&\text{for}\qquad
\g=\oa_{2n+1}\\[0.3em]
\la_n(a)\ts \si_n(a)&=\la_{n+1}(a)\ts \si_n(a+2),\quad\qquad&&\text{for}\qquad
\g=\spa_{2n}\\[0.3em]
\la_{n-1}(a)\ts \si_n(a)&=\la_{n+1}(a)\ts \si_n(a+1),\quad\qquad&&\text{for}\qquad
\g=\oa_{2n}.
\end{alignedat}
\eeq

For every $i\in\{1,\dots,n\}$ define a linear operator $\wt S_i:\Lc\to \wt\Lc_i$
satisfying the Leibniz rule \eqref{lr}.
For $i=1,\dots,n-1$ set
\beql{defsibcd}
\wt S_i:\la_j(a)\mapsto\begin{cases} \la_i(a)\ts\si_i(a)\qquad&\text{for}\quad j=i\\
-\la_{i+1}(a)\ts\si_i(a+1)\qquad&\text{for}\quad j=i+1\\
-\la_{i'}(a)\ts\si_i(a+\ka-i+1)\qquad&\text{for}\quad j=i'\\
\la_{(i+1)'}(a)\ts\si_i(a+\ka-i)\qquad&\text{for}\quad j=(i+1)'\\
0\qquad&\text{for}\quad j\ne i,i',i+1,(i+1)'.
\end{cases}
\eeq
The action of $\wt S_n$ depends on the type and is given as follows.

\paragraph{Case $\g=\oa_{2n+1}${\rm :}} $\wt S_n:\la_j(a)\mapsto 0$ if $j<n$ or $j>n'$ and

\ben
\bal
\la_n(a)&\mapsto\la_n(a)\ts\big(\si_n(a)+\si_n(a-1/2)\big)\\[0.3em]
\la_{n+1}(a)&\mapsto\la_{n+1}(a)\ts\big(\si_n(a-1/2)-\si_n(a+1/2)\big)\\[0.3em]
\la_{n'}(a)&\mapsto -\la_{n'}(a)\ts\big(\si_n(a)+\si_n(a+1/2)\big).
\eal
\een

\paragraph{Case $\g=\spa_{2n}${\rm :}} $\wt S_n:\la_j(a)\mapsto 0$ if $j<n$ or $j>n'$ and

\ben
\bal
\la_n(a)&\mapsto\la_n(a)\ts \si_n(a)\\[0.3em]
\la_{n'}(a)&\mapsto -\la_{n'}(a)\ts \si_n(a+2).
\eal
\een

\paragraph{Case $\g=\oa_{2n}${\rm :}} $\wt S_n:\la_j(a)\mapsto 0$ if $j<n-1$ or $j>(n-1)'$ and

\ben
\bal
\la_{n-1}(a)&\mapsto\la_{n-1}(a)\ts \si_n(a)\\[0.3em]
\la_n(a)&\mapsto\la_n(a)\ts \si_n(a)\\[0.3em]
\la_{n'}(a)&\mapsto -\la_{n'}(a)\ts \si_n(a+1)\\[0.3em]
\la_{(n-1)'}(a)&\mapsto -\la_{(n-1)'}(a)\ts \si_n(a+1).
\eal
\een

The relations \eqref{larel} are easily seen to be preserved by the action
of the $\wt S_i$ so that
the operators on $\Lc$ are well-defined.
The {\em $i$-th screening operator}
\ben
S_i:\Lc\to \Lc_i
\een
is now defined as the composition of $\wt S_i$ and the projection $\wt\Lc_i\to \Lc_i$.

Due to \cite[Theorem~5.1]{fm:cq}, we can define the subalgebra $\Rep\Y(\g)$
of Yangian characters in $\Lc$ as the intersection of kernels of the screening
operators:
\ben
\Rep\Y(\g)=\bigcap_{i=1}^{n} \ts\ker S_i.
\een

\bre\label{rem:notlbcd}
The variables $\la_i(a)$ and $\si_i(a)$
are related to the corresponding elements used in \cite{fr:qc}
as follows: $\la_i(a)=\La_{i,q^{4a}}$ for $\g=\oa_{2n+1}$
and $\la_i(a)=\La_{i,q^{2a}}$ for $\g=\spa_{2n}$ and $\g=\oa_{2n}$.
Moreover, for $\g=\oa_{2n+1}$
\ben
\si_i(a)=S_{i,q^{4a+2i-2}}\qquad\text{for}\quad i=1,\dots,n-1;\qquad
\si_n(a)=S_{n,q^{4a+2n-1}},
\een
while for $\g=\spa_{2n}$ we have
\ben
\si_i(a)=S_{i,q^{2a+i-1}}\qquad\text{for}\quad i=1,\dots,n;
\een
the latter relations with $i=1,\dots,n-1$ hold for $\g=\oa_{2n}$ as well,
but $\si_n(a)=S_{n,q^{2a+n-2}}$; cf. \cite{nn:pt}.
Note also that relations \eqref{larel} were obtained in
\cite[Prop.~5.2 and 5.14]{amr:rp} as the conditions for the highest weight
representations of the Yangian $\Y(\g)$ to be nontrivial, whereas \eqref{rellasibcd}
and \eqref{quo} are consistent with the conditions on the representation to
be finite-dimensional; cf. \cite[Theorem 5.16]{amr:rp}.
\qed
\ere

We follow \cite[Sec.~8.1]{f:lc} again
to define the classical $\Wc$-algebra $\Wc(\g)$;
see also \cite{mm:yc} and \cite{mr:cw}.
We will regard $\U(\wh\h_-)$
as the algebra of polynomials in the variables $F_{ii}[r]$
with $i=1,\dots,n$ and $r<0$. The screening operators
\ben
V_i:\U(\wh\h_-)\to \U(\wh\h_-),\qquad i=1,\dots,n
\een
are defined as follows. For $i=1,\dots,n-1$ set
\ben
V_i=\sum_{r=0}^{\infty} V_{i\ts [r]}\ts
\Bigg(\frac{\di}{\di F_{ii}[-r-1]}-\frac{\di}{\di F_{i+1\ts i+1}[-r-1]}\Bigg),
\een
where the coefficients $V_{i\ts [r]}$ are found from the expansion
of a formal generating function in a variable $z$,
\ben
\sum_{r=0}^{\infty} V_{i\ts [r]}\ts z^r=\exp\ts\sum_{m=1}^{\infty}
\frac{F_{ii}[-m]-F_{i+1\ts i+1}[-m]}{m}\ts z^m.
\een
For the action of $V_n$ we have the following formulas.

\paragraph{Case $\g=\oa_{2n+1}${\rm :}}

\ben
V_n=\sum_{r=0}^{\infty} V_{n\ts [r]}\ts
\frac{\di}{\di F_n[-r-1]},
\een
where
\ben
\sum_{r=0}^{\infty} V_{n\ts [r]}\ts z^r=\exp\ts\sum_{m=1}^{\infty}
\frac{F_n[-m]}{m}\ts z^m.
\een

\paragraph{Case $\g=\spa_{2n}${\rm :}}

\ben
V_n=\sum_{r=0}^{\infty} V_{n\ts [r]}\ts
\frac{\di}{\di F_n[-r-1]},
\een
where
\ben
\sum_{r=0}^{\infty} V_{n\ts [r]}\ts z^r=\exp\ts\sum_{m=1}^{\infty}
\frac{2\tss F_n[-m]}{m}\ts z^m.
\een

\paragraph{Case $\g=\oa_{2n}${\rm :}}

\ben
V_n=\sum_{r=0}^{\infty} V_{n\ts [r]}\ts
\Bigg(\frac{\di}{\di F_{n-1}[-r-1]}+\frac{\di}{\di F_n[-r-1]}\Bigg)
\een
where
\ben
\sum_{r=0}^{\infty} V_{n\ts [r]}\ts z^r=\exp\ts\sum_{m=1}^{\infty}
\frac{F_{n-1}[-m]+ F_n[-m]}{m}\ts z^m.
\een

The {\em classical $\Wc$-algebra} $\Wc(\g)$ is a subalgebra
of $\U(\wh\h_-)$ defined as the intersection of kernels of the screening
operators:
\ben
\Wc(\g)=\bigcap_{i=1}^{n} \ts\ker V_i.
\een

Now construct a map $\gr:\Rep\Y(\g)\to \Wc({}^L\g)$ and describe its properties.
First, embed $\Lc$ into the algebra of formal power series $\CC[[\la^{(r)}_i]]$
in variables $\la^{(r)}_i$ with $i=1,\dots,N$ and $r=0,1,\dots$ by using
\eqref{embedla} and taking the quotient by the corresponding relations
\eqref{larel}. Introduce new variables $\mu^{(r)}_i$
by \eqref{lamu} for $i=1,\dots,n$ and
define their degrees by $\deg\mu_i^{(r)}=-r-1$.
Given $A\in\Lc$, consider the corresponding element $\CC[[\mu^{(r)}_i]]$
and take its homogeneous component $\overline A$ of the maximum degree.
This component is a polynomial in the variables $\mu^{(r)}_i$
and so we have a map
\beql{grbcd}
\gr: \Lc\to \CC[\mu^{(r)}_i],\qquad A\mapsto\overline A.
\eeq
Note its property \eqref{grab}. We will identify
$\U(\wh\h_-)$ with the algebra of polynomials $\CC[\mu^{(r)}_i]$ via the isomorphism
$F_{ii}[-r-1]\mapsto \mu^{(r)}_i/r!$.

\bpr\label{prop:grtypebcd}
The image of the restriction of the map \eqref{grbcd} to the subalgebra
of characters $\Rep\Y(\g)$ is contained in $\Wc({}^L\g)$ and so it defines a map
\ben
\gr:\Rep\Y(\g)\to \Wc({}^L\g).
\een
Moreover, any homogeneous element
of $\Wc({}^L\g)$ is contained in the image of $\gr$.
\epr

\bpf
The proof is quite similar to that of Proposition~\ref{prop:grtypea}
so we only point out the changes to be made. Introduce variables $\si^{(r)}_i$ by
\eqref{embedsi}
and set $\deg\si^{(r)}_i=-r-1$. Define operators
$S^{\tss\circ}_i$ for $i=1,\dots,n$ as in \eqref{scirc}, where the quotient
is now taken by the respective relations \eqref{rellasibcd}
and \eqref{quo} written in terms of the $\mu^{(r)}_i$ with $a$ understood
as a variable. Since relations \eqref{rellasibcd}
are identical to \eqref{rellasi}, the argument for the operators
$S^{\tss\circ}_i$ with $i=1,\dots,n-1$ follows the same steps as for
type $A$. To complete the proof for the operator $S^{\tss\circ}_n$,
consider the three cases separately.

\paragraph{Case $\g=\oa_{2n+1}${\rm .}} As with \eqref{barsi},
the corresponding operator $\overline S_n$ is
now given by
\beql{barsib}
\overline S_n:\mu^{(r)}_j\mapsto\begin{cases} 2\tss\si^{(r)}_n\qquad&\text{for}\quad j=n\\
0\qquad&\text{for}\quad j\ne n.
\end{cases}
\eeq
Note that
\ben
\la_{n+1}(a)=\frac{\la_1(a+n-1)\la_2(a+n-2)\dots \la_n(a)}{\la_1(a+n-1/2)\la_2(a+n-3/2)\dots \la_n(a+1/2)}
\een
which is easy to derive from \eqref{larel}. Now use \eqref{quo} and write $\la_i(a)=1+\mu_i(a)$
for $i=1,\dots,n$ to get the corresponding analogue of \eqref{mua}.
As a result, we get the equation
\ben
\si'_n(z)=2\tss\mu_n(z)\ts \si_n(z)
\een
so that for the images under
$\overline S_n$ we have
\ben
\overline S_n:\mu_n(z)\mapsto 2\tss\exp 2\int\mu_n(z)\ts dz,
\een
and $\overline S_n:\mu_j(z)\mapsto 0$ for $j\ne n$. This coincides
with the action of the operator $2\tss V_n$ associated with $\spa_{2n}$ on the series
\ben
\mu_n(z)=\sum_{r=0}^{\infty} F_{nn}[-r-1]\ts z^r.
\een
Hence, if an element $A\in\Lc$ is annihilated
by all operators $S_i$, then its image $\overline A$ under the map \eqref{gr} is
annihilated by all operators $V_i$
associated with $\spa_{2n}$ which is Langlands dual to $\oa_{2n+1}$.

\paragraph{Case $\g=\spa_{2n}${\rm .}} Similar to \eqref{barsib},
we have
\beql{barsic}
\overline S_n:\mu^{(r)}_j\mapsto\begin{cases} \si^{(r)}_n\qquad&\text{for}\quad j=n\\
0\qquad&\text{for}\quad j\ne n.
\end{cases}
\eeq
Relations \eqref{larel} now imply
\ben
\la_{n+1}(a)=\frac{\la_1(a+n)\la_2(a+n-1)\dots \la_{n-1}(a+2)}{\la_1(a+n+1)\la_2(a+n)\dots \la_n(a+2)}.
\een
Write $\la_i(a)=1+\mu_i(a)$
for $i=1,\dots,n$ and use \eqref{quo} to get the equation
\ben
\si'_n(z)=\mu_n(z)\ts \si_n(z).
\een
Hence, for the images under
$\overline S_n$ we have
\ben
\overline S_n:\mu_n(z)\mapsto \exp \int\mu_n(z)\ts dz,
\een
and $\overline S_n:\mu_j(z)\mapsto 0$ for $j\ne n$. This coincides
with the action of the operator $V_n$ associated with $\oa_{2n+1}$ on the series
\ben
\mu_n(z)=\sum_{r=0}^{\infty} F_{nn}[-r-1]\ts z^r.
\een
Therefore, if an element $A\in\Lc$ is annihilated
by all operators $S_i$, then its image $\overline A$ under the map \eqref{gr} is
annihilated by all operators $V_i$
associated with $\oa_{2n+1}$ which is Langlands dual to $\spa_{2n}$.

\paragraph{Case $\g=\oa_{2n}${\rm .}} Similar to \eqref{barsi},
we have
\beql{barsid}
\overline S_n:\mu^{(r)}_j\mapsto\begin{cases} \si^{(r)}_n\qquad&\text{for}\quad j=n-1,n\\
0\qquad&\text{for}\quad j\ne n-1,n.
\end{cases}
\eeq
We derive from \eqref{larel} that
\ben
\la_{n+1}(a)=\frac{\la_1(a+n-2)\la_2(a+n-3)\dots \la_{n-1}(a)}{\la_1(a+n-1)\la_2(a+n-2)\dots \la_n(a)}.
\een
Write $\la_i(a)=1+\mu_i(a)$
for $i=1,\dots,n$ and use \eqref{quo} to get the equation
\ben
\si'_n(z)=\big(\mu_{n-1}(z)+\mu_n(z)\big)\ts \si_n(z).
\een
Hence, for the images under
$\overline S_n$ we have
\ben
\overline S_n:\mu_{n-1}(z)\mapsto \exp\int\big(\mu_{n-1}(z)+\mu_n(z)\big)\ts dz,\qquad
\mu_n(z)\mapsto \exp\int\big(\mu_{n-1}(z)+\mu_n(z)\big)\ts dz,
\een
and $\overline S_n:\mu_j(z)\mapsto 0$ for $j\ne n-1,n$. This coincides
with the action of the operator $V_n$ associated with $\oa_{2n}$ on the series
\ben
\mu_i(z)=\sum_{r=0}^{\infty} F_{ii}[-r-1]\ts z^r,\qquad i=n-1,n.
\een
Thus, if an element $A\in\Lc$ is annihilated
by all operators $S_i$, then its image $\overline A$ under the map \eqref{gr} is
annihilated by all operators $V_i$
associated with $\oa_{2n}$ which is Langlands self-dual.

The last part of the proposition follows from \cite{mm:yc}, where generators of the
classical $\Wc$-algebra
were obtained as images of the Yangian characters
under the map $\gr$.
\epf


\end{document}